\documentclass[conference]{IEEEtran}

\ifCLASSINFOpdf
\else
\fi
\usepackage{graphics} 
\usepackage{epsfig} 
\usepackage{amsmath} 
\usepackage{amssymb}  
\usepackage[latin1]{inputenc}
\usepackage{amsfonts}
\usepackage{subfigure}

\newtheorem{lemme}{Lemma}
\newtheorem{corollaire}{Corollary}
\newtheorem{prop}{Proposition}

\begin{document}

\title{Mathematical properties of a semi-classical signal analysis method: noisy signal case}

\author{\IEEEauthorblockN{Da-Yan Liu}
\IEEEauthorblockA{Mathematical and \\ Computer Sciences and \\ Engineering Division,\\
King Abdullah University \\of Science and Technology (KAUST), KSA\\
Email: dayan.liu@kaust.edu.sa}
\and
\IEEEauthorblockN{Taous-Meriem Laleg-Kirati}
\IEEEauthorblockA{Mathematical and \\ Computer Sciences and \\ Engineering Division,\\
King Abdullah University \\of Science and Technology (KAUST), KSA\\
Email: taousmeriem.laleg@kaust.edu.sa}}

\maketitle

\begin{abstract}
Recently, a new signal analysis method based on a semi-classical approach  has been proposed \cite{meriem}. The main idea in this method  is to interpret a signal as a potential of a Schrodinger operator and then to use the discrete spectrum of this operator to analyze the signal. In this paper, we are interested  in  a mathematical analysis  of this method in  discrete case considering noisy signals.
\end{abstract}

\IEEEpeerreviewmaketitle

\section{INTRODUCTION}
Recently, a new signal analysis method based on a semi-classical approach  has been proposed \cite{meriem}. We refer to this method SCSA for  \textbf{S}emi-\textbf{C}lassical \textbf{S}ignal \textbf{A}nalysis. The main idea in the SCSA is to interpret a signal as a potential of a Schrodinger operator depending on a semi-classical parameter \cite{Dimassi}. It is well-known that if the potential is in the Faddeev class \cite{Faddeev}, then it can be expressed using a sum of the squared eigenfunctions associated to the negative eigenvalues characterizing the discrete spectrum  of the Schrodinger operator and an integral involving  the continuous spectrum. Similarly to the other standard approximation methods, by truncating the expression of the potential, the sum part is taken as an estimate of the potential. The proposed estimate depends on the semi-classical parameter. It has been  shown that by reducing this parameter the estimation of the signal by the SCSA can be improved \cite{meriem}.


Promising results have been obtained when applying the SCSA to  arterial blood pressure. More than a satisfactory estimation of the pressure signals, this method introduced new spectral parameters that seem  to contain important physiological information \cite{meriem07,meriem10}. Moreover, a recent study has shown that the SCSA parameters  could  be useful in the  estimation of  some physical parameters related to turbomachinery features \cite{Fadi}. It has been also confirmed by some tests that the SCSA method is robust with respect to corrupting noises. Hence, it  can be considered as a filter.
However, the   mathematical analysis of the SCSA in  discrete noisy case  has not been considered yet, comparing to other signal analysis methods  (see, e.g., \cite{ans,Liu2008,Liu2009,Liu2011b,Liu2011c}),

We propose in this paper to study the mathematical properties of the SCSA in  discrete case considering noisy signals.  We propose also to study how to  choose an appropriate semi-classical parameter to analyze a signal.
In Section \ref{section1}, we  recall the methodology and some mathematical  properties of the SCSA method in the continuous case. Then, we study  this method in  discrete case in Section \ref{section2}.
In Section \ref{section3},  the SCSA method is studied in  discrete  noisy case. An a-posteriori error bound of the noise error contribution  is given. In Section \ref{section4}, we show how to choose an appropriate parameter for the SCSA method in  discrete noisy case. Moreover, we can see its efficiency  and stability with corrupting noises. Finally, we give some conclusions and
perspectives for our future work in Section \ref{section4}.

\section{PRELIMINARY}\label{section1}

In this section, we are going to recall the methodology and some mathematical  properties of the SCSA  \cite{meriem}.

\subsection{Methodology}

Let us consider the following Schrodinger operator
 \begin{equation}\label{Eq_schrodinger_operator}
    H_h(y):=-h^2 \frac{d^2}{d x^2} - y,
 \end{equation}
where $h \in \mathbb{R}^*_+$ is a semi-classical parameter \cite{Dimassi}, $y \in \mathcal{L}^1_1(\mathbb{R}):= \left\{V | \int_{-\infty}^{+\infty} \left| V(x)\right|(1+|x|)\,dx < \infty  \right\}$, and
$\mathcal{L}^1_1(\mathbb{R})$ is called the Faddeev class \cite{Faddeev}. Then, it is well known that $y$ can be reconstructed as follows \cite{Kaup1978}
\begin{equation}\label{Eq_Deift_Trubowitz}
\begin{split}
    & \quad \quad  y(x) =  4h \sum_{n=1}^{N_h} \kappa_{n_h} \psi_{n_h}^2(x)- \\
    & \lim_{b\rightarrow +\infty}\frac{1}{b} \int_0^b \left( \frac{2i}{\pi}h \int_{-a}^{+a} k R_{r(l)h}(k)\,f^2_{\pm h}(k,x)\,dk \right)\,da, \ a.e.,
\end{split}
\end{equation}
where $- \kappa_{n_h}^2$ are the negative eigenvalues of $H_h(y)$ with  $\kappa_{1_h}>\kappa_{2_h}>\cdots>\kappa_{N_h}>0$, $N_h$ denotes the number of the negative eigenvalues, and $\psi_{n_h} \in \mathcal{H}^2(\mathbb{R})$ ($\mathcal{H}^2(\mathbb{R})$ being the Sobolev space of order $2$)  are the associated  $\mathcal{L}^2$-normalized eigenfunctions such that
 \begin{equation}\label{}
    H_h(y)\,\psi_{n_h}=-\kappa_{n_h}^2\,\psi_{n_h}.
 \end{equation}
Moreover,
$R_{r(l)h}(\cdot)$ is the reflection coefficient and $f_{h\pm}(\cdot,\cdot)$ denote the Jost solutions defined as the unique solutions of the Schrodinger integral equation at $\pm \infty$ respectively.

Let us mention that the eigenpairs  $\left\{\kappa_{n_h},\psi_{n_h}\right\}$ for $n_h=1,\cdots, N_h$ can be numerically calculated, but   $R_{r(l)h}$ and $f_{h\pm}$ usually can not be calculated. Hence, similarly to the other standard approximation methods, by truncating (\ref{Eq_Deift_Trubowitz}) we propose to take the sum part as an estimate of $y$. Moreover, since this sum part $4h \displaystyle\sum_{n=1}^{N_h} \kappa_{n_h} \psi_{n_h}^2$ is a positive definite function on $\mathbb{R}$, it is necessary to assume $y$
to be positive definite too. Then, we give the following proposition.

\begin{prop}\cite{meriem} \label{Prop_scsa}
Let $y$  be a positive definite  function belonging to $\mathcal{L}^1_1(\mathbb{R})$, then it can be estimated by
\begin{equation}\label{Eq_Estimation0}
    y_h(x):= 4h \sum_{n=1}^{N_h} \kappa_{n_h} \psi_{n_h}^2(x),
\end{equation}
where $h>0$, $- \kappa_{n_h}^2$ are the negative eigenvalues of the Schrodinger operator $H_h(y)$ defined by (\ref{Eq_schrodinger_operator}) with  $\kappa_{1_h}>\kappa_{2_h}>\cdots>\kappa_{N_h}>0$, $N_h$ denotes the number of the negative eigenvalues, and $\psi_{n_h}$  are the associated  $\mathcal{L}^2$-normalized eigenfunctions.  \hfill$\Box$
\end{prop}

Consequently, the proposed estimation only depends on the parameter $h$ \cite{meriem}.

\subsection{Some properties of the SCSA method}

Now, it is natural to consider the convergence of this proposed method which is shown in the following proposition.
\begin{prop}\cite{meriem} \label{Prop_accuracy}
Let  $y$ be a real valued function satisfying  the following condition
\begin{equation}\label{Eq_space_B}
\begin{split}
& \quad \quad \quad y \in  \mathcal{B}:= \\ &\left\{V \in \mathcal{L}^1_1(\mathbb{R})|  \ \forall\, x \in \mathbb{R}, \ y(x) \geq 0, \ \frac{\partial^m y}{\partial x^m} \in  \mathcal{L}^1(\mathbb{R}),\  m=1,2 \right\}.
    \end{split}
\end{equation}
Then, we have
$\lim_{h \rightarrow 0}\left\|y_h-y\right\|_{\mathcal{L}^1(\mathbb{R})}=0$.

Moreover, the number $N_h$ of the negative eigenvalues of $H_h(y)$ is a decreasing function of $h$.  \hfill$\Box$
\end{prop}

According to the previous proposition, we can improve our estimation by reducing the value of $h$. Moreover, it was shown in \cite{meriem} that if there exists an $h$ such that $\frac{y}{h^2}$ is a reflectionless potential of the Schrodinger operator $H_h(y)$, then $y_h$ is an exact representation of $y$.

Let us recall that the study of the Schrodinger operator in the case where $h \rightarrow 0$ is referred to the semi-classical
analysis \cite{Dimassi}. Consequently, we call the proposed signal estimation method  \emph{\textbf{S}emi-\textbf{C}lassical \textbf{S}ignal \textbf{A}nalysis (SCSA) method}.




\section{The SCSA method in  discrete case}\label{section2}

In order to apply the SCSA method, we need to solve numerically the following Schrodinger eigenvalue problem
 \begin{equation}\label{Eq_schrodinger_eigenvalue_problem}
    H_h(y)\,\psi_{h}=\lambda_{h}\,\psi_{h}, \ \ \psi_{h} \in \mathcal{H}^2(\mathbb{R}).
 \end{equation}
Hence, we will  study  the SCSA method in  discrete case in  this section.

We assume that $y^{d}$ is a discrete  signal defined on an interval $I=[a,b] \subset\mathbb{R}$,
which can be considered as the restriction of a function $y$ satisfying the condition given in (\ref{Eq_space_B}), such that
 \begin{equation}\label{Eq_siganl_discrete}
y^{d}(x_j) = y(x_j),
 \end{equation}
 where $x_j=a + (j-1)\, \Delta x$ for $j=1,\cdots, M$, are equidistant points with
the distance between two consecutive points $\Delta x= \frac{b-a}{M-1}$. We denote $y^{d}_j$ as the value  $y^{d}(x_j)$  for $j=1,\cdots, M$.

The discretization of the Schrodinger eigenvalue problem given in (\ref{Eq_schrodinger_eigenvalue_problem}) leads to the following eigenvalue matrix problem
\begin{equation}\label{Eq_matrix1}
    A^d_h  \hat{\Psi}_{h} = \hat{\lambda}_{h} \hat{\Psi}_{h},
\end{equation}
where $\hat{\Psi}_{h}=\left[\hat{\psi}_{h,1},\hat{\psi}_{h,2},\cdots,\hat{\psi}_{h,M-1},\hat{\psi}_{h,M} \right]^T$,
\begin{equation} \label{Eq_matrix11}
A^d_h=-h^2 D_2-\text{diag}(Y^{d}),
\end{equation}
$\text{diag}(Y^{d})$ is a diagonal matrix whose elements are $y^{d}_j$ for $j=1,\cdots, M$, and
$D_2$ is a second order differentiation matrix given by a discretization method for differential equations \cite{Boyd,Trefethen}, which is independent of $h$.

Let us consider the negative eigenvalues $-\hat{\kappa}^2_{n_h}$ of $A^d_h$ with $\hat{\kappa}_{1_h} \geq \cdots \geq \hat{\kappa}_{\hat{N}_h}>0$, where $\hat{N}_h$ is the number of the negative eigenvalues of $A^d_h$ with $0 \leq \hat{N}_h \leq M$. We denote their associated eigenvectors by $\hat{\Psi}_{n_h}=\left[\hat{\psi}_{n_h,1},\cdots,\hat{\psi}_{n_h,M} \right]^T$ for $n_h=1,\cdots, \hat{N}_h$. Moreover, we assume that
\begin{equation}\label{Eq_norm1}
\Delta x\sum_{j=1}^{M} \hat{\psi}^2_{n_h,j}=1, \quad i.e. \quad \left\|\hat{\Psi}_{n_h}\right\|_2 = \frac{1}{\sqrt{\Delta x}}.
\end{equation}
Thus, according to Proposition \ref{Prop_scsa} we can construct an estimation of $y$ by the SCSA method in this discrete case as follows
\begin{equation}\label{Eq_Estimation1}
y^{d}_h(x_j):= 4h \sum_{n=1}^{\hat{N}_h} \hat{\kappa}_{n_h} \hat{\psi}_{n_h,j}^2,
\end{equation}
for $j=1,\cdots, M$.
 By writing the following equality
\begin{equation}
\begin{split}
y^{d}-y^{d}_h = \left(y^{d}- y_h\right) + \left(y_h- y^{d}_h\right),
\end{split}
\end{equation}
we can see that the estimation error in this discrete case for $y^{d}_h$  can be divided into two sources:
\begin{enumerate}
  \item the estimation error for $y_h$:  $y^{d} - y_h$, which corresponds to the truncated integral part in (\ref{Eq_Deift_Trubowitz}),
  \item the discrete numerical error for  $y^{d}_h$: $y_h-y^{d}_h$, which is produced by the discretization of the Schrodinger eigenvalue problem.
\end{enumerate}
It is shown in Proposition \ref{Prop_accuracy} that the number of the negative eigenvalues of the Schrodinger operator is decreasing with respect to $h$. This property is generalized to the discrete case in the following proposition.

\begin{prop} \label{Prop_propeties2}
We assume that the matrix $A^d_h=-h^2 D_2 -\text{diag}(Y^{d})$ defined  in (\ref{Eq_matrix11}) satisfies the following conditions:
\begin{description}
  \item[{(C1):}]  $D_2$ is  symmetric,
  \item[{(C2):}]  $D_2$ is  negative definite,
  \item[{(C3):}]  the number of zeros in the diagonal of $\text{diag}(Y^{d})$ is equal to $\hat{n}$ with $0 \leq \hat{n} \leq M-1$.
\end{description}
Let us denote  the number  of the negative eigenvalues of $A^d_h$ by $\hat{N}_h$, then we have
\begin{equation}\label{Eq_boundNh}
\forall\, h>0,\ 0\leq \hat{N}_h \leq M - \hat{n},
\end{equation}
where $M$ is the size of the matrix $A^d_h$. Moreover, we have
\begin{align}
 \hat{N}_{h} &= 0, \ \text{ for } h > \sqrt{\frac{\hat{y}^d_{M}}{d_{M}}},\label{Eq_infNh}\\
 \hat{N}_{h} &= M - \hat{n}, \ \text{ for } 0< h < \sqrt{\frac{\hat{y}^d_{\hat{n}+1}}{d_{1}}}\label{Eq_supNh},
\end{align}
where $\hat{y}^d_{\hat{n}+1}$ (resp. $\hat{y}^d_{M}$) is the smallest (resp. largest) strictly  positive element in the diagonal of $\text{diag}(Y^{d})$,  and $d_{1}$ (resp. $d_{M}$) is the largest (resp. smallest) eigenvalue of $D_2$.
   \hfill$\Box$
\end{prop}

\noindent{\textbf{Proof.}} Since the matrices $D_2$ and $\text{diag}(Y^{d})$ defined in (\ref{Eq_matrix11}) both are symmetric, $A^d_h$ is also
symmetric. The eigenvalues of $A^d_h$, $-D_2$ and $-\text{diag}(Y^{d})$ can be denoted and ordered respectively as follows
\begin{equation*}\label{}
\hat{\lambda}_{h,M} \leq \cdots \leq \hat{\lambda}_{h,1}, \ d_{M} \leq \cdots \leq d_{1},\ -\hat{y}^d_{M} \leq \cdots \leq -\hat{y}^d_{1}.
\end{equation*}
Then, according to Weyl's theorem (see \cite{Stewart})
we get
\begin{equation}\label{Eq_ordered_eigenvalues}
h^2 d_{M}-\hat{y}^d_{j} \leq  \hat{\lambda}_{h,j} \leq h^2 d_{1}-\hat{y}^d_{j},
\end{equation}
for $j=1,\cdots,M$. If $\hat{n}=0$, then (\ref{Eq_boundNh}) can be directly obtained.  By using (C2) and (C3) we obtain
\begin{equation}\label{Eq_ordered_eigenvalues1}
\forall\, h>0, \ 0 < h^2 d_{M} \leq  \hat{\lambda}_{h,j},
\end{equation}
for $j=1,\cdots,\hat{n}$, and
\begin{equation}\label{Eq_ordered_eigenvalues2}
\forall\, h>0, \ \hat{\lambda}_{h,j} \leq h^2 d_{1}-\hat{y}^d_{j},
\end{equation}
for $j=\hat{n}+1,\cdots,M$.
Hence, (\ref{Eq_boundNh}) can  be deduced from (\ref{Eq_ordered_eigenvalues1}).
Moreover, if $0 \leq \hat{n} \leq M-1$, then by using  (\ref{Eq_ordered_eigenvalues2}) we obtain that for any  $h < \sqrt{\frac{\hat{y}^d_{\hat{n}+1}}{d_{1}}}$,  $\hat{\lambda}_{h,j}<0$, for $j=\hat{n}+1,\cdots,M$. It yields $\hat{N}_{h} \geq M - \hat{n}$. Then, (\ref{Eq_infNh}) can be obtained by using (\ref{Eq_boundNh}).
 Finally, this proof can be completed by solving  the following inequality
 \begin{equation}\label{Eq_ordered_eigenvalues}
0< h^2 d_{M}-\hat{y}^d_{M} \leq  \hat{\lambda}_{h,M}.
\end{equation}
 \hfill$\Box$

\begin{corollaire} \label{}
We assume that the matrix $A^d_h$ defined  in (\ref{Eq_matrix11}) satisfies the  conditions (C1)-(C3) given in Proposition \ref{Prop_propeties2}.
Then,  $\forall\, h>0$, $\exists\, h'$ with $0 < h' < h$, such that
\begin{equation}\label{Eq_inegalityNh}
\hat{N}_{h'} \geq \hat{N}_{h},
\end{equation}
where $\hat{N}_{h'}$ and $\hat{N}_{h}$ denote  the number  of the negative eigenvalues of $A^d_{h'}$  and $A^d_h$ respectively.    \hfill$\Box$
\end{corollaire}

\noindent{\textbf{Proof.}} Let $h>0$, by using (\ref{Eq_boundNh}) we have $0 \leq \hat{N}_{h} \leq M - \hat{n}$.  Then, by using (\ref{Eq_supNh}) we get
 \begin{equation}
\forall\, 0< h'<\min\left(h, \sqrt{\frac{\hat{y}^d_{\hat{n}+1}}{d_{1}}}\right), \
\hat{N}_{h'} \geq \hat{N}_{h}.
\end{equation}
Thus, this proof is completed.
\hfill$\Box$

\section{ERROR ANALYSIS IN   NOISY CASE}\label{section3}

In this section, we are going to consider the SCSA method in  discrete  noisy case. Moreover, an a-posteriori error bound of the noise error contribution for the SCSA method will be given.

\subsection{The SCSA method in  discrete noisy case}From now on, we assume that
\begin{equation}\label{Eq_signal_noisy}
y^{\varpi}=y^{d}+\varpi
\end{equation}
is a  noisy observation of the discrete signal $y^{d}$ defined in (\ref{Eq_siganl_discrete}),
where the corrupting noise $\varpi$ is an identically distributed sequence of random variables with an expected value $\mu$  and a variance  $\sigma^2$ ($\sigma \in \mathbb{R}_+$).
We denote $y^{\varpi}_j$ and $\varpi_j$ as the values  $y^{\varpi}(x_j)$ and $\varpi(x_j)$ respectively,  for $j=1,\cdots, M$.

In order to apply the SCSA method in this discrete noisy case, by substituting $A^d_h$ in (\ref{Eq_matrix1}) by $A^{\varpi}_h$ we need to solve the following eigenvalue matrix problem
\begin{equation}\label{Eq_matrix2}
    A^{\varpi}_h  \tilde{\Psi}_{n} = \tilde{\lambda}_{h} \tilde{\Psi}_{h},
\end{equation}
where
\begin{equation} \label{Eq_matrix22}
A^{\varpi}_h=-h^2 D_2-\text{diag}(Y^{\varpi}),
\end{equation}
$\text{diag}(Y^{\varpi})$ is a diagonal matrix whose  elements are $y^{\varpi}_j$ for $j=1,\cdots, M$. Hence, according to (\ref{Eq_matrix11}) we obtain
\begin{equation} \label{Eq_matrix_relation}
A^{\varpi}_h= A^d_h -\text{diag}(W),
\end{equation}
$\text{diag}(W)$ is a diagonal matrix whose  elements  are $\varpi_j$ for $j=1,\cdots, M$.

Let us denote the negative eigenvalues of $A^{\varpi}_h$ by $-\tilde{\kappa}^2_{n_h}$ with  $\tilde{\kappa}_{1_h} \geq \cdots \geq \tilde{\kappa}_{\tilde{N}_h}>0$ for $n_h=1,\cdots, \tilde{N}_h$, where $\tilde{N}_h$ is the number of the negative eigenvalues of $A^{\varpi}_h$. Similar to (\ref{Eq_norm1}),  their associated eigenvectors are denoted by  $\tilde{\Psi}_{n_h}=\left[\tilde{\psi}_{n_h,1},\cdots,\tilde{\psi}_{n_h,M} \right]^T$ with
$\left\|\tilde{\Psi}_{n_h}\right\|_2 = \frac{1}{\sqrt{\Delta x}}$.
Thus, according to Proposition \ref{Prop_scsa}, an estimation of $y$  can be given by the SCSA method  in this discrete noisy case as follows
\begin{equation}\label{Eq_Estimation2}
    y^{\varpi}_h(x_j):= 4h \sum_{n=1}^{\tilde{N}_h} \tilde{\kappa}_{n_h} \tilde{\psi}_{n_h,j}^2,
\end{equation}
for $j=1,\cdots, M$.
By writing the following equality
\begin{equation}
\begin{split}
y^{d}-y^{\varpi}_h = \left(y^{d} - y_h\right) + \left(y_h- y^{d}_h\right) + \left(y^{d}_h - y^{\varpi}_h\right),
\end{split}
\end{equation}
we can see that the total estimation error in this discrete noisy case for $y^{\varpi}_h$ can be divided into three parts:
\begin{enumerate}
  \item the truncated error for $y_h$: $y^{d} - y_h$,
  \item the discrete numerical error for $y^{d}_h$:  $y_h - y^{d}_h$,
  \item the noise error contribution for $y^{\varpi}_h$: $y^{d}_h - y^{\varpi}_h$.
\end{enumerate}

Proposition \ref{Prop_propeties3} illustrates a property on the number of the negative eigenvalues  in
the discrete noisy case.

\begin{prop} \label{Prop_propeties3}
We assume that the matrix $A^{\varpi}_h=-h^2 D_2 -\text{diag}(Y^{\varpi})$ defined  in (\ref{Eq_matrix22}) satisfies the conditions (C1)-(C2) and the following condition
\begin{description}
  \item[{(C4):}]  the number of the positive elements in the diagonal of $-\text{diag}(Y^{\varpi})$ is equal to $\tilde{n}$ with $0 \leq \tilde{n} \leq M-1$.
\end{description}
Let us denote  the number  of the negative eigenvalues of $A^{\varpi}_h$ by $\tilde{N}_h$, then we have
\begin{equation}\label{Eq_boundNh2}
\forall\, h>0,\ 0\leq \tilde{N}_h \leq M - \tilde{n},
\end{equation}
\begin{align}
\tilde{N}_{h} &= 0, \text{ for } \  h > \sqrt{\frac{\tilde{y}^{\varpi}_{M}}{d_{M}}},\label{Eq_infNh2}\\
\tilde{N}_{h} &= M - \tilde{n}, \text{ for } \  0< h < \sqrt{\frac{\tilde{y}^{\varpi}_{\tilde{n}+1}}{d_{1}}},\label{Eq_supNh2}
\end{align}
where  $\tilde{y}^{\varpi}_{\tilde{n}+1}$ (resp. $\tilde{y}^{\varpi}_{M}$) is the smallest (resp. largest) strictly  positive element in the diagonal of $\text{diag}(Y^{\varpi})$, and $d_{1}$ (resp. $d_{M}$) is the largest (resp. smallest) eigenvalue of $D_2$.

Moreover, we have $\forall\, h>0$, $\exists\, h'$ with $0 < h' < h$, such that
$\tilde{N}_{h'} \geq \tilde{N}_{h},$
where $\tilde{N}_{h'}$ and $\tilde{N}_{h}$ denote  the number  of the negative eigenvalues of $A^{\varpi}_{h'}$  and $A^{\varpi}_h$ respectively.   \hfill$\Box$
\end{prop}

\noindent{\textbf{Proof.}} This proof can be completed in a  similar way to the one of Proposition \ref{Prop_propeties2}.
\hfill$\Box$

\subsection{Analysis of the noise error contribution} \label{subsection2}

In this subsection, we are going to study the noise error contribution in the SCSA method by providing an a-posteriori error bound in the following proposition.

\begin{prop} \label{Prop_error}
Let $y^{\varpi}$ be a discrete noisy signal defined as in (\ref{Eq_signal_noisy}),   $y^{\varpi}_h$ and $y^{d}_h$ be the estimations of $y$ given by (\ref{Eq_Estimation1}) and (\ref{Eq_Estimation2}) respectively. Moreover, we assume that
\begin{description}
  \item[{(C1):}] the matrix $D_2$ given in (\ref{Eq_matrix11}) and (\ref{Eq_matrix22}) is symmetric,
  \item[{(C5):}] the numbers of the negative eigenvalues of $A^d_h$ and $A^{\varpi}_h$ are equal, $i.e.$   $\hat{N}_{h}=\tilde{N}_{h}$,
  \item[{(C6):}]
$\tilde{\kappa}^2_{n_h} < 2\, \hat{\kappa}^2_{n_h}$,  for $n_h=1,\cdots,\tilde{N}_h$.
\end{description}
Then, an a-posteriori error bound for the noise error contribution in $y^{\varpi}_h$ can be given as follows
 \begin{equation}
\begin{split}
 & \left\|y^{\varpi}_h - y^{d}_h\right\|_2
\stackrel{p}{<}  \frac{4h}{\sqrt{\Delta x}} \sum_{n=1}^{\tilde{N}_h} \left( 2 \tilde{\kappa}_{n_h} + \frac{B^{\gamma}_{\mu,\sigma}}{\sqrt{2}\,\tilde{\kappa}_{n_h}} \right),
\end{split}
\end{equation}
where $B^{\gamma}_{\mu,\sigma}=\max\left(|\mu-
\gamma \sigma|,|\mu+
\gamma \sigma|\right)$ with $\mu$ and $\sigma^2$ being the expected value and the variance of $\varpi_j$ respectively for $j=1,\cdots, M$, and  $c \stackrel{p}{<} d$ means that the probability for $c$ to be smaller than $d$ is $p$ with $p=1- \frac{1}{\gamma^2}$ and $\gamma \in \mathbb{R}^*_+$. \hfill$\Box$
\end{prop}

In order to prove this proposition, we need the following lemma.
\begin{lemme} \label{Lem}
By giving the conditions $\text{(C1)}$, $\text{(C5)}$,  and $\text{(C6)}$,  we have
\begin{equation}
\begin{split}
\left| \tilde{\kappa}_{n_h}- \hat{\kappa}_{n_h} \right|
& \stackrel{p}{<} \frac{B^{\gamma}_{\mu,\sigma}}{\sqrt{2}\,\tilde{\kappa}_{n_h}},
\end{split}
\end{equation}
for $n_h=1,\cdots,\tilde{N}_h$, where  $B^{\gamma}_{\mu,\sigma}=\max\left(|\mu-
\gamma \sigma|,|\mu+
\gamma \sigma|\right)$.

 \hfill$\Box$
\end{lemme}

\noindent{\textbf{Proof.}} According to the condition  $\text{(C1)}$, the matrixes $A^d_h$ and $A^{\varpi}_h$ are both symmetric. Moreover, since the matrix $\text{diag}(W)$ is diagonal, its eigenvalues are $\varpi_j$ for $j=1,\cdots,M$.
Hence, by using Weyl's theorem (see \cite{Stewart}) and (\ref{Eq_matrix_relation}) with the condition  $\text{(C5)}$, we obtain
\begin{equation} \label{Eq_difference_kappa2}
\left| \tilde{\kappa}^2_{n_h}- \hat{\kappa}^2_{n_h} \right| \leq \max_{1 \leq j \leq M} \left| \varpi_j \right|,
\end{equation}
for $n_h=1,\cdots,\tilde{N}_h$.
By using the Bienaym\'{e}-Chebyshev inequality we get that for any real number $\gamma>0$,
\begin{equation}
\text{Pr}\left(\left|\varpi_j-\mu\right| < \gamma \sigma\right)
> 1-\frac{1}{\gamma^2},
\end{equation}
$i.e.$ the probability for $\varpi_j$ to be
within the interval  $\left]\mu-
\gamma  \sigma \,, \mu+
\gamma  \sigma\right[$  is higher than
 $1- \frac{1}{\gamma^2}$. Consequently, by denoting $B^{\gamma}_{\mu,\sigma}=\max\left(|\mu-
\gamma \sigma|,|\mu+
\gamma \sigma|\right)$ and using (\ref{Eq_difference_kappa2}) we obtain
\begin{equation} \label{Eq_bond_noise}
\left| \tilde{\kappa}^2_{n_h}- \hat{\kappa}^2_{n_h} \right| \leq \max_{1 \leq j \leq M} \left| \varpi_j \right| \stackrel{p}{<} B^{\gamma}_{\mu,\sigma},
\end{equation}
for $n_h=1,\cdots,\tilde{N}_h$. According to the condition $(\text{C6})$ we obtain
 \begin{equation}
 0 < \frac{\tilde{\kappa}^2_{n_h}}{2} < \min\left(\tilde{\kappa}^2_{n_h},\hat{\kappa}^2_{n_h}\right).
 \end{equation}
   Then, the utilization of the mean value theorem gives us
\begin{equation} \label{Eq_difference_kappa}
\begin{split}
 &\left| \tilde{\kappa}_{n_h}- \hat{\kappa}_{n_h} \right| \leq
 \frac{1}{2} \frac{ \left| \tilde{\kappa}^2_{n_h}- \hat{\kappa}^2_{n_h} \right|}{ \left(  \frac{\tilde{\kappa}^2_{n_h}}{2}  \right)^{\frac{1}{2}}}.
 \end{split}
\end{equation}
Then, this proof can be completed by using (\ref{Eq_bond_noise}).
\hfill$\Box$

\noindent{\textbf{Proof of Proposition \ref{Prop_error}.}}
By using the condition $\text{(C5)}$, (\ref{Eq_Estimation1}) and (\ref{Eq_Estimation2}), we get
\begin{equation}\label{Eq_difference_y}
\begin{split}
y^{\varpi}_h(x_j) - y^{d}_h(x_j) & = 4h \sum_{n=1}^{\tilde{N}_h} \left(\tilde{\kappa}_{n_h} \tilde{\psi}_{n_h,j}^2 -  \hat{\kappa}_{n_h} \hat{\psi}_{n_h,j}^2\right),
\end{split}
\end{equation}
for $j=1,\cdots,M$. Let us denote $E_{n_h}=[e_{n_h,1}, \cdots, e_{n_h,M}]^T$ with
\begin{equation} \label{Eq_difference_kappa}
e_{n_h,j} =  \tilde{\kappa}_{n_h} \tilde{\psi}_{n_h,j}^2 -  \hat{\kappa}_{n_h} \hat{\psi}_{n_h,j}^2.
\end{equation}
Then, we obtain
\begin{equation}
\begin{split}
e_{n_h,j}
&= \tilde{\kappa}_{n_h} \tilde{\psi}_{n_h,j}^2 -\tilde{\kappa}_{n_h}\hat{\psi}_{n_h,j}^2
+ \tilde{\kappa}_{n_h} \hat{\psi}_{n_h,j}^2 -  \hat{\kappa}_{n_h} \hat{\psi}_{n_h,j}^2\\
&= \tilde{\kappa}_{n_h} \left(\tilde{\psi}_{n_h,j}^2-\hat{\psi}_{n_h,j}^2  \right) + \left( \tilde{\kappa}_{n_h}- \hat{\kappa}_{n_h} \right) \hat{\psi}_{n_h,j}^2.
\end{split}
\end{equation}
By calculating  the norms $\left\|\hat{\Psi}_{n_h}\right\|_2$ and $\left\|\tilde{\Psi}_{n_h}\right\|_2$,  $\left\|E_{n_h}\right\|_2$ can be  bounded as follows
\begin{equation}\label{Eq_E}
\left\|E_{n_h}\right\|_2 \leq  \frac{2 \tilde{\kappa}_{n_h} + \left| \tilde{\kappa}_{n_h}- \hat{\kappa}_{n_h} \right|}{\sqrt{\Delta x} }.
\end{equation}
Hence, by using (\ref{Eq_difference_y}), (\ref{Eq_difference_kappa}) and (\ref{Eq_E}), we get
\begin{equation}
\begin{split}
\left\|y^{\varpi}_h - y^{d}_h\right\|_2 & \leq 4h \sum_{n=1}^{\tilde{N}_h} \left\|E_{n_h}\right\|_2\\
&  \leq \frac{4h}{\sqrt{\Delta x}} \sum_{n=1}^{\tilde{N}_h} \left( 2 \tilde{\kappa}_{n_h}+  \left| \tilde{\kappa}_{n_h}- \hat{\kappa}_{n_h} \right|\right).
\end{split}
\end{equation}
Finally, the proof can be completed by using Lemma \ref{Lem}.
\hfill$\Box$

The convergence of the SCSA method in the continuous noise-free case is shown in Proposition \ref{Prop_accuracy}, where
the estimation proposed by the SCSA method can be improved by reducing the value of $h$. Now, it is interesting to study
the  efficiency  of this method in  discrete noisy case. Especially, we need to know  the influence of $h$ on the noise error contribution in the SCSA method.
A natural idea is to study the influence of the parameter $h$ on
the noise error bound given in Proposition \ref{Prop_error} so as to deduce the one on the noise error contribution.
On one hand, because of the term $h$ this error bound seems to be increasing with respect to $h$. On the other hand, it is shown in Proposition \ref{Prop_propeties3}
that the number of the negative eigenvalues $\tilde{N}_h$ can be  decreasing with respect to $h$.  Consequently, it is impossible to  intuitively know the influence of $h$. However, it can be studied numerically as shown in  the next section.


\section{NUMERICAL RESULTS} \label{section4}

In this section, by taking a numerical example we are going to show how to choose an appropriate value of $h$ for the SCSA method in  discrete noisy case. Moreover, we can see its efficiency and stability with corrupting noises.

\begin{figure}[h!]
\centering {\includegraphics[scale=0.48]{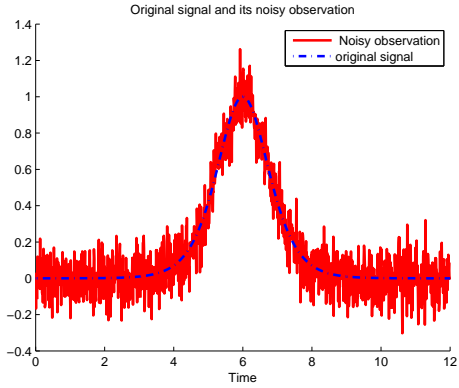}}\caption{The discrete signal $y^{d}$ and the discrete noisy signal $y^{\varpi}$ with $SNR=11\text{dB}$.}%
\label{Fig_signal}%
\end{figure}

We assume that $y^{\varpi}$ is the discrete noisy
observation of $y$ defined in (\ref{Eq_signal_noisy}),  where
\begin{equation}\label{}
y^{d}(x_j)=y(x_j)= \text{sech}^2(x_j-6)
\end{equation}
with $x_j \in I=[0,12]$ and $\Delta x=10^{-2}$.
Hence, $M$ is equal to $1201$. Moreover, we assume that the noise
$\varpi$ is simulated from a zero-mean white Gaussian $iid$ sequence.  The variance $\sigma^2$ of
$\varpi$ is adjusted in such a way that  the signal-to-noise ratio is equal to
$11\text{dB}$.  We can see the original signal $y$ and its noisy observation $y^{\varpi}$ in Figure \ref{Fig_signal}.
Let us recall that this sech-squared  function is well known in the quantum physics theory as the Pöschl-Teller potential of the Schrodinger operator \cite{Olafsson}.

In order to estimate $y$ by using  (\ref{Eq_Estimation2}),
we propose to use a Fourier pseudo-spectral method \cite{Boyd} to solve numerically the Schrodinger eigenvalue problem defined in (\ref{Eq_schrodinger_eigenvalue_problem}). Thus,
the second order differentiation matrix $D_2$ is  given as follows \cite{Trefethen}:
 If $M$ is even, then
  \begin{equation*}\label{}
    D_2(k,j)=\frac{\Delta}{\left(\Delta x\right)^2} \left\{
                                                      \begin{array}{ll}
                                                        -\frac{\pi^2}{3\Delta^2}-\frac{1}{6}, &  \text{for } k = j,  \\
                                                        -(-1)^{k-j}\frac{1}{2}\frac{1}{\sin^{2}\left(\frac{(k-j)\Delta}{2}\right)}, &  \text{for } k \neq j.
                                                      \end{array}
                                                    \right.
\end{equation*}
If $M$ is odd, then
  \begin{equation*}\label{}
    D_2(k,j)=\frac{\Delta}{\left(\Delta x\right)^2} \left\{
                                                      \begin{array}{ll}
                                                        -\frac{\pi^2}{3\Delta^2}-\frac{1}{12}, &  \text{for } k = j,  \\
                                                        -(-1)^{k-j}\frac{1}{2}\frac{\cot \left(\frac{(k-j)\Delta}{2}\right)}{\sin \left(\frac{(k-j)\Delta}{2}\right)}, &  \text{for } k \neq j.
                                                      \end{array}
                                                    \right.
\end{equation*}
with $\Delta=\frac{2\pi}{M}$. Let us mention  that $D_2$
is symmetric and negative  definite. 
Then, we use the Matlab routine
$eig$ to solve the eigenvalue matrix problems defined in (\ref{Eq_matrix1}) and (\ref{Eq_matrix2}).

It is shown in Subsection \ref{subsection2} that the total estimation error for $y^{\varpi}_h$ comes from three parts. However, since the
estimation $y_h$ can not be calculated in the discrete case, we only consider the estimation error in the discrete noise-free case and the noise error contribution
\begin{equation}
y^{d}-y^{\varpi}_h=\left(y^{d} - y^{d}_h\right)+ \left(y^{d}_h - y^{\varpi}_h\right).
\end{equation}

\begin{figure}[h!]
\centering {\includegraphics[scale=0.48]{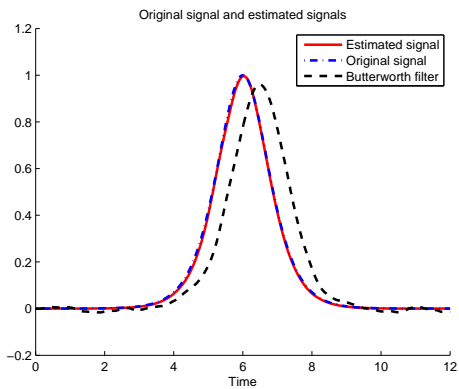}}\caption{The discrete signal $y^{d}$, the signal ${}_{f}y^{\varpi}$  filtered by Butterworth filter, and the estimate $y^{\varpi}_h$ obtained by the SCSA method with $h=0.4$.}%
\label{Fig_estimation}%
\end{figure}

\begin{figure}[h!]
 \centering
 \subfigure[The $l^2$ norms of different errors.]
 {\includegraphics[scale=0.295]{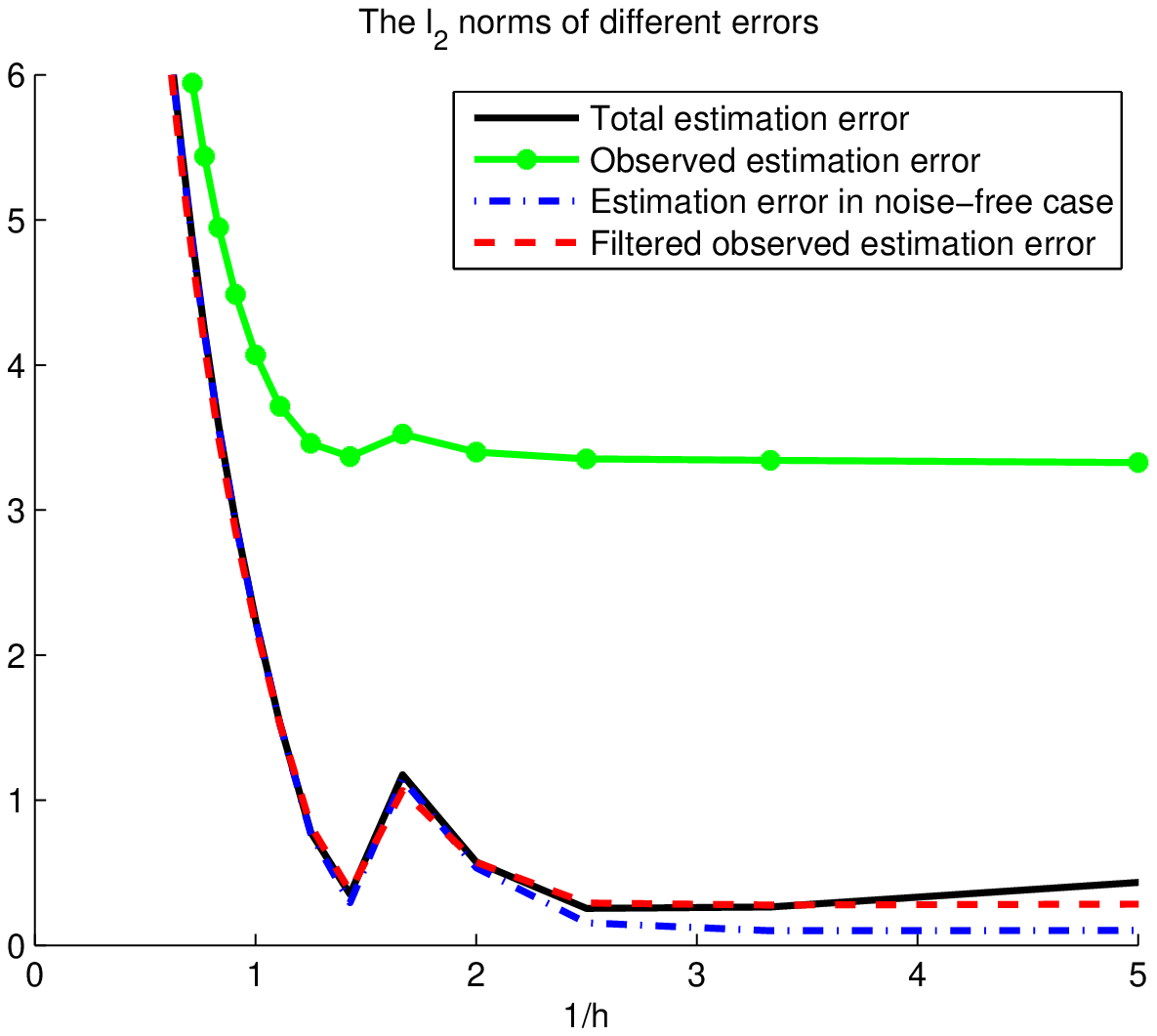}\label{Fig_error1}}
 \subfigure[The noise error contribution.]
 {\includegraphics[scale=0.295]{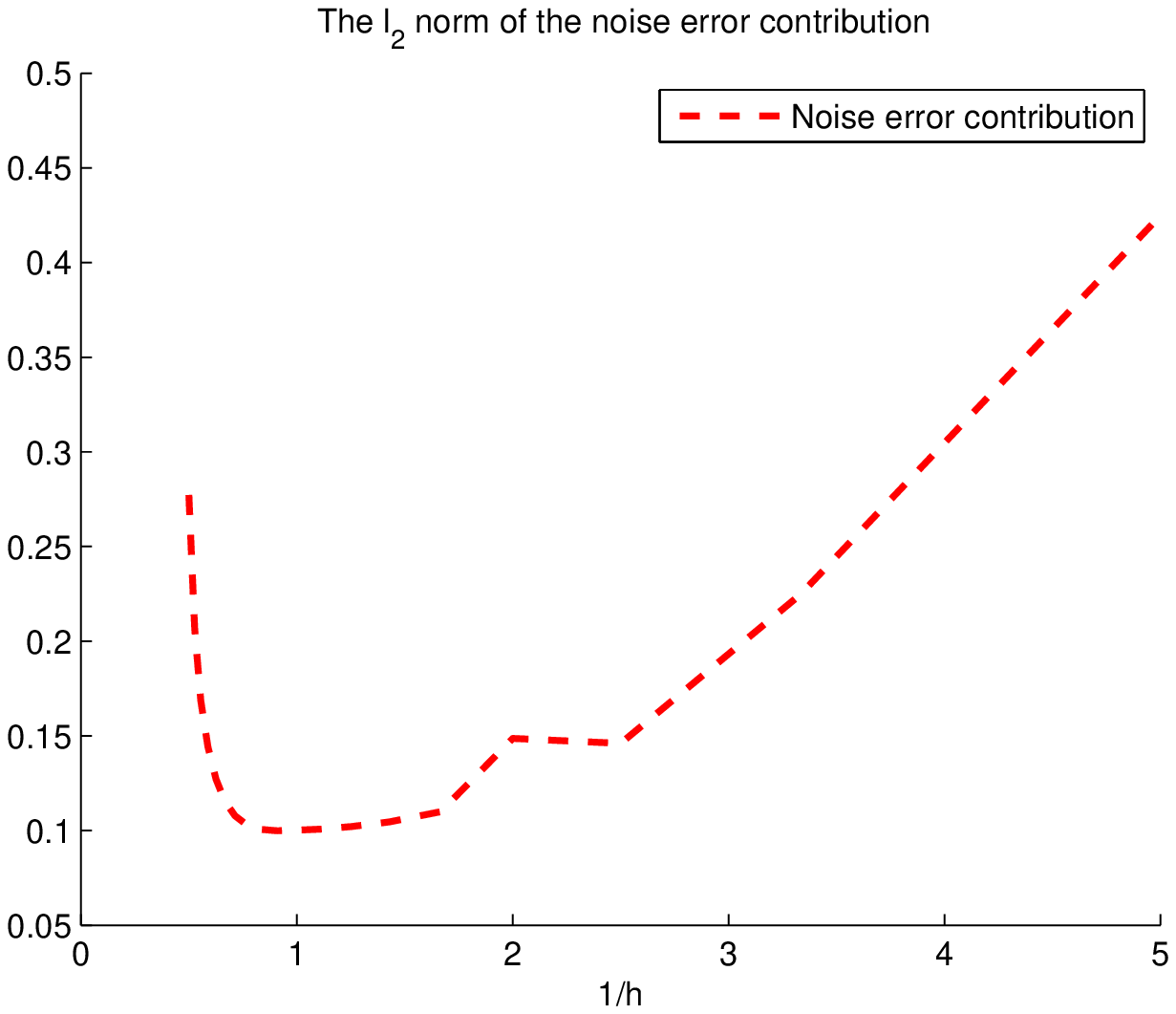}\label{Fig_error2}}
   \subfigure[The noise error bound]
 {\includegraphics[scale=0.295]{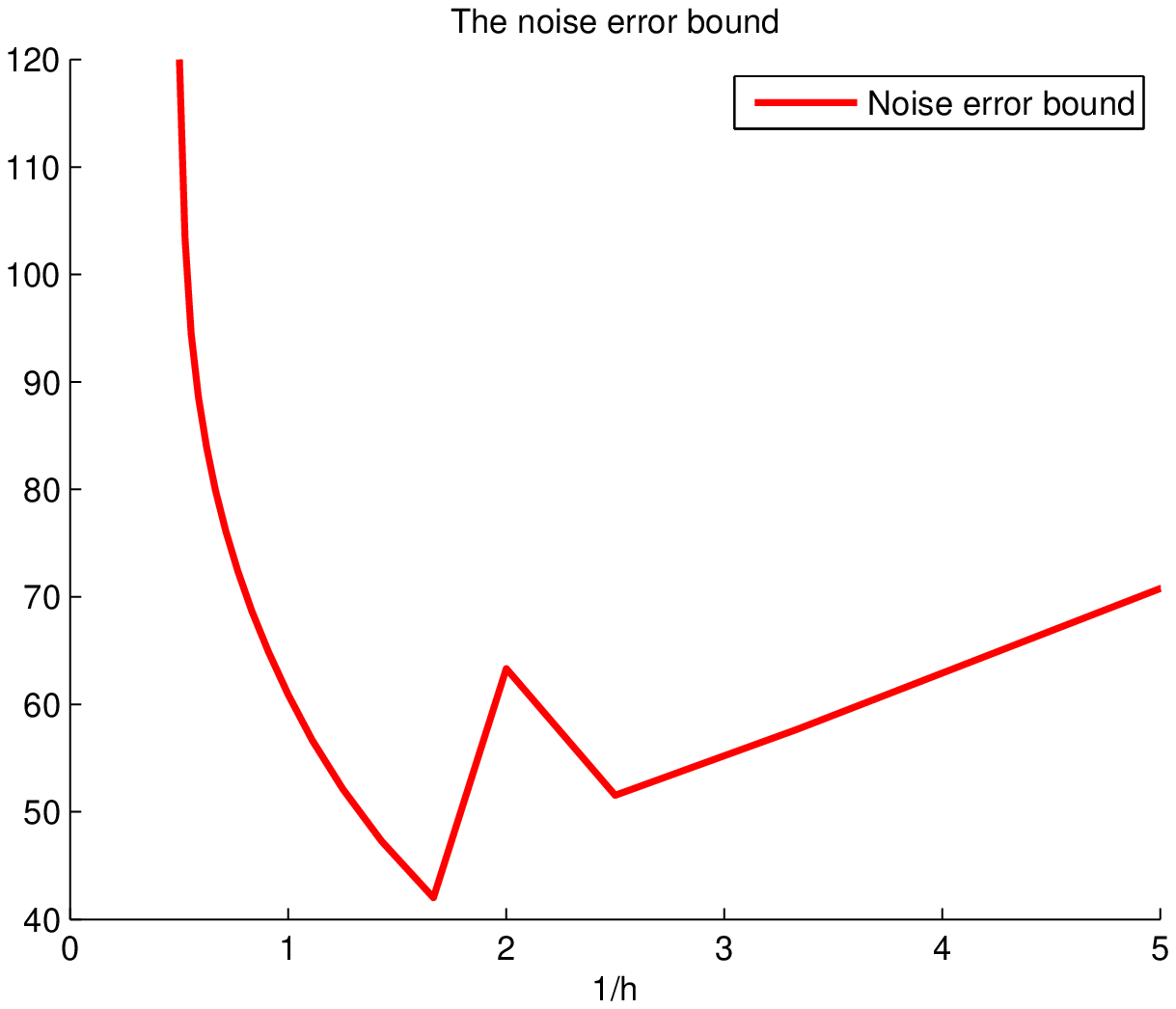}\label{Fig_errorbound}}
  \subfigure[The number of negatives eigenvalues.]
 {\includegraphics[scale=0.295]{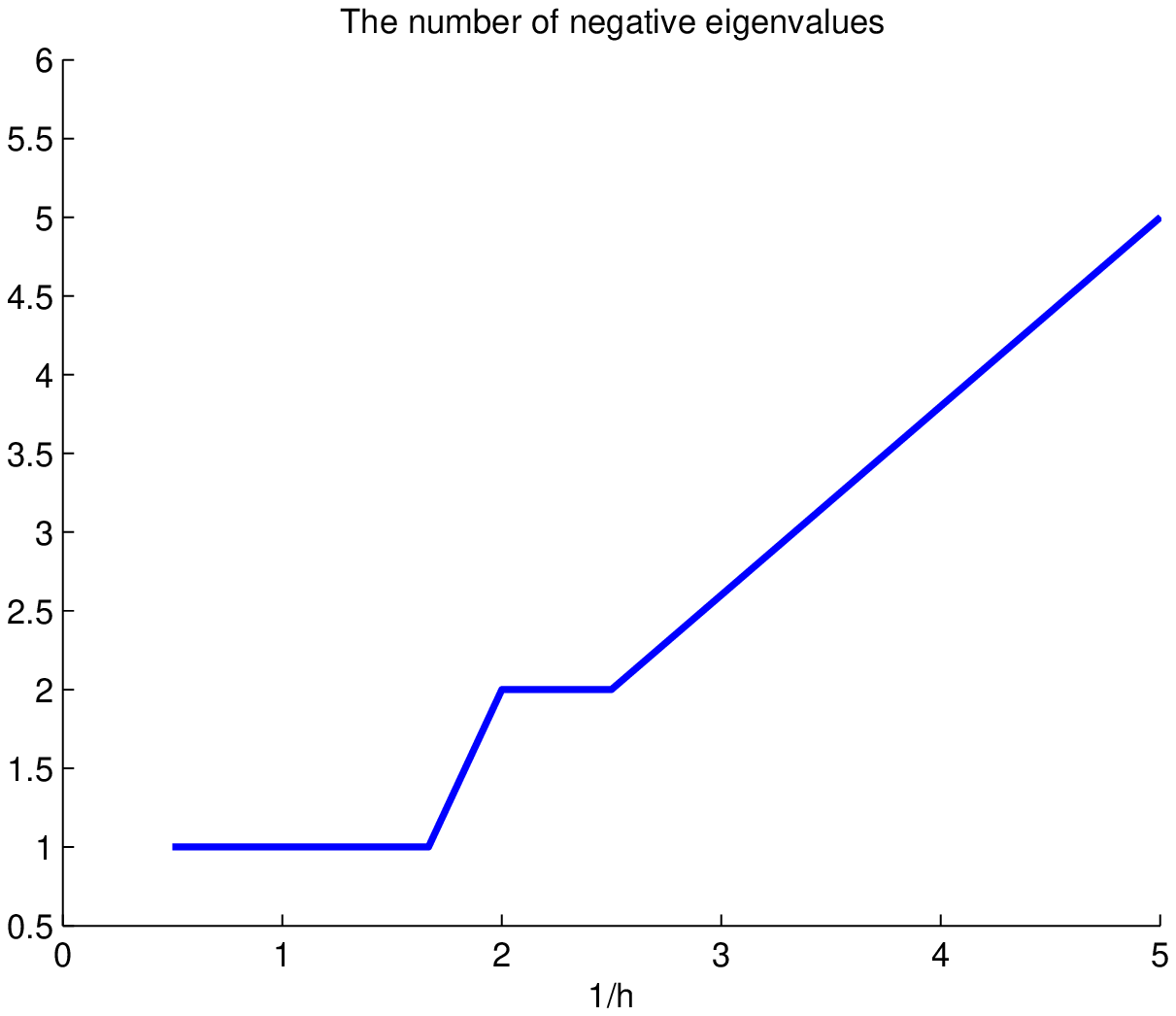}\label{Fig_nb}}
\caption{Results of different variations obtained  for $h=0.2,0.3,\cdots,1.9,2$.} \label{}
\end{figure}

In order to see the influence of $h$ on the total estimation error,  we show the variations of
$\left\|y^{d}-y^{\varpi}_h\right\|_2$ in Figure \ref{Fig_error1}, which is represented by the black solid line. We can see that $\left\|y^{d}-y^{\varpi}_h\right\|_2$ has a minimum at $h=0.4$ and a local minimum at  $h=0.7$. Thus, we can take the optimal value  $h=0.4$ for $y^{\varpi}_h$ so as to produce a minimal total  estimation error.
The estimation obtained  by using $y^{\varpi}_h$ with $h=0.4$ is given  in Figure \ref{Fig_estimation}. Consequently, we can see that the SCSA method is accurate and robust with a corrupting noise. Hence, it can be considered as a filter for noisy signals without delays.

The previous analysis is based on the knowledge of $y^{d}$ which is usually unknown in the practice work. If we use $\left\|y^{\varpi}-y^{\varpi}_h\right\|_2$ to study the influence of $h$ on the total estimation error, then we generally can not find an optimal value of $h$. This can be explained by the green solid line in Figure  \ref{Fig_error1}, which corresponds to the different values  of $\left\|y^{\varpi}-y^{\varpi}_h\right\|_2$.  In order to solve this problem, we propose to use a second-order Butterworth filter \cite{Bianchi} which is given as follows
\begin{equation}
H(s)=\frac{w^2_c}{s^2+2 w_c s +w_c^2},
\end{equation}
where the cutoff frequency $w_c$ is set to $w_c=0.01$.

$H(\cdot)$ is a classical low-pass filter which can be used to attenuate the corrupting Gaussian noise in $y^{\varpi}$. The filtered signal is represented by the black dotted line  in Figure \ref{Fig_estimation}.  Since this filter produces a delay to $y^{\varpi}$, we also apply it to $y^{\varpi}_h$. The filtered signals are denoted by ${}_{f}y^{\varpi}$ and ${}_{f}y^{\varpi}_h$ respectively. Then,  we use the influence of $h$ on $\left\|{}_{f}y^{\varpi}-{}_{f}y^{\varpi}_h\right\|_2$  to deduce the one on $\left\|y^{d}-y^{d}_h\right\|_2$ which corresponds to the estimation error in the discrete noise-free case. We can see in Figure  \ref{Fig_error1} the relation between the variations of $\left\|{}_{f}y^{\varpi}-{}_{f}y^{\varpi}_h\right\|_2$ and $\left\|y^{d}-y^{d}_h\right\|_2$ which are represented by
the red dotted line and the blue dash-dotted line respectively. Hence, we can observe that they have the same variation with the same local minimum and the same local maximum.

Now, we study the influence of $h$ on the noise error contribution. The variations of
$\left\|y^{d}_h-y^{\varpi}_h\right\|_2$ is shown in Figure \ref{Fig_error2}. Moreover, we can verify that
the number of the negative eigenvalues $\hat{N}_h$ and $\tilde{N}_h$ are equal. We can see their variation with respect to $h$ in Figure \ref{Fig_nb}.
Consequently, according to Figure \ref{Fig_error2} and Figure \ref{Fig_nb} we can observe  that when $\hat{N}_h$ is equal to $1$ the noise error contribution is increasing with respect to $h$. However, when the value of $\hat{N}_h$ increases, the noise error contribution is decreasing with respect to $h$. Consequently, a small value of $h$ can produce a large noise error contribution.

We are going to use Proposition \ref{Prop_error} to deduce the variation of the noise error contribution with respect to $h$.
Since the noise is assumed to be a  zero-mean white Gaussian $iid$ sequence, by using the well known three-sigma rule we obtain that
\begin{equation} \label{}
 \max_{1 \leq j \leq M} \left| \varpi_j \right| \stackrel{99.7\%}{<} 3\sigma.
\end{equation}
The variation of the noise error bounds given in Proposition \ref{Prop_error} is shown in Figure \ref{Fig_errorbound}. Although this error bound is not sharp, its variation is similar to the one of the noise error contribution shown in Figure \ref{Fig_error2}, where there are local minimums  at $h=0.6$ and  $h=0.4$, and a local maximum at $h=0.5$.

Finally, by combined with the variation of $\left\|{}_{f}y^{\varpi}-{}_{f}y^{\varpi}_h\right\|_2$ shown in  Figure  \ref{Fig_error1}, we can
choose $h=0.4$ or $h=0.7$ in our estimation so as to minimize the total estimation error.

\section{conclusion} \label{section5}

In this paper, the SCSA method recently introduced for signal analysis is studied in  discrete noisy case.
Some mathematical properties of the negative eigenvalues of a Schrodinger operator are given in  discrete noise-free case and  discrete noisy case respectively. An a-posteriori error bound of the noise error contribution is given which is based on the expected value and the variance of a corrupting noise. By taking a numerical example, we show the influence of the semi-classical parameter on different sources of errors in the SCSA method. Moreover, it is shown that the SCSA method is accurate and robust against corrupting noises. Hence, it can be considered as a filter without involved delays.  Finally, we study how to choose an appropriate semi-classical parameter
without knowing the original signal. The comparison to other signal analysis methods like Fourier transform or the wavelets will be done  in a future work. Moreover, the SCSA method will be extended for time derivatives filtering, which is still an open problem, such that this method can be useful in more applications in signal processing and automatic control.



\end{document}